\documentclass[11pt]{article}

\usepackage{latexsym,mathrsfs}
\usepackage{amsmath,amssymb}
\usepackage{amsthm,enumerate,verbatim}
\usepackage{amsfonts}
\usepackage{graphicx}
\usepackage{algorithm}
\usepackage{algorithmic}
\usepackage{url}

\newtheorem{lemma}{Lemma}

\newtheorem{theorem}{Theorem}

\newtheorem{definition}{Definition}
\newtheorem*{definition*}{Definition}

\numberwithin{equation}{section}
\numberwithin{table}{section}
\numberwithin{figure}{section}
\def \R{{\mathbb R}}
\def \C{{\mathbb C}}
\def\real{\mathop{\mathrm{Re}}}
\def\imag{\mathop{\mathrm{Im}}}

\newcommand {\mat}  [1] {\left[\begin{array}{#1}}
\newcommand {\rix}      {\end{array}\right]}

\DeclareMathOperator{\argmin}{argmin}

\def\real{\mathop{\mathrm{Re}}}
\def\imag{\mathop{\mathrm{Im}}}
\newcommand{\eproof}{\space
    {\ \vbox{\hrule\hbox{\vrule height1.3ex\hskip0.8ex\vrule}\hrule}}\par}

\title{On approximating the nearest $\Omega$-stable matrix }
\date{}

\author{Neelam Choudhary\thanks{
 School of Engineering and Applied Sciences, Department of Mathematics, Bennett University, Greater Noida-201310, Uttar Pradesh, India; \texttt{neelam.choudhary@bennett.edu.in. }}
 \qquad Nicolas Gillis\thanks{
 Department of Mathematics and Operational Research,
Facult\'e Polytechnique, Universit\'e de Mons, Rue de Houdain~9, 7000 Mons, Belgium; \texttt{nicolas.gillis@umons.ac.be}. N. Gillis
acknowledges the support of the ERC (starting grant n$^\text{o}$ 679515) and F.R.S.-FNRS (incentive grant for scientific research n$^\text{o}$ F.4501.16).
}
 \qquad Punit Sharma\thanks{Department of Mathematics, Indian Institute of Technology Delhi, Hauz Khas, New Delhi-110016, India; \texttt{punit.sharma@maths.iitd.ac.in}.
P. Sharma acknowledges the support of the DST-Inspire Faculty Award (MI01807-G) by Government of India and Institute SEED Grant (NPN5R) by IIT Delhi.
 }
}

\begin{document}

\maketitle

\begin{abstract}
In this paper, we consider the problem of approximating a given matrix with a matrix whose eigenvalues lie in some specific region $\Omega$ of the complex plane. More precisely, we consider three types of regions and their intersections: conic sectors, vertical strips and disks. We refer to this problem as the nearest $\Omega$-stable matrix problem. This includes as special cases the stable matrices for continuous and discrete time linear time-invariant systems. In order to achieve this goal, we parametrize this problem using dissipative Hamiltonian matrices and linear matrix inequalities. This leads to a reformulation of the problem with a convex feasible set. By applying a block coordinate descent method on this reformulation, we are able to compute solutions to the approximation problem, which is illustrated on some examples.
\end{abstract}


\section{Introduction}

Let us consider the following linear time-invariant (LTI) systems of the form
\[
\dot{x}(t)=Ax(t)+Bu(t),
\]
where $A\in \mathbb R^{n,n}$, $B\in \mathbb R^{n,m}$, $x$ is the state vector, $u$ is the input vector.
Let also $\Omega$ be a subset of the complex plane.
An LTI system is called $\Omega$-stable if all eigenvalues of $A$ lie inside $\Omega$.
The most well-known example is when $\Omega$ is the open left-half of the complex plane which characterizes stability of continuous LTI systems. For discrete LTI systems of the form
$x(l+1) = Ax(l)$ for $l \in \mathbb N$, stability requires $\Omega$ to be the unit disk.
Another motivation for enforcing eigenvalues in specific regions of the complex plane is that
the transient response is related to the location of its eigenvalues~\cite{Kuo87,Acker93}.
For example, the step response of
a system with eigenvalues $\lambda=-\tau \omega_n\pm i\omega_d$ is fully characterized in terms of the undamped natural frequency $\omega_n=|\lambda|$,
the damping ratio $\tau$, and the damped natural frequency $\omega_d$~\cite{Kuo87}.
By constraining $\lambda$ to lie in a prescribed region,
specific bounds can be put on these quantities to ensure a satisfactory transient response; see also~\cite{trefethen2005spectra} and the references therein.
In this paper, we focus on three regions of the complex plane and their intersections, namely:
\begin{itemize}
\item Conic sector: the conic sector region of parameters $a,\,\theta \in \R$ with $0 \leq \theta \leq \pi/2$, denoted by $\Omega_C(a,\theta)$, is defined as
\[
\Omega_C(a,\theta)
:=
\left\{ x+iy \in \C \ \big| \ \sin(\theta) (x-a) < \cos(\theta) y < -\sin(\theta) (x-a),\, x \leq a \right\}.
\]

\item Vertical strip: the vertical strip region of parameters $ h < k$, denoted by $\Omega_V(h,k)$, is defined as
\[
\Omega_V(h,k)
:= \left\{ x+iy \in \C \ \big| \ -k < x < -h \right\}.
\]
Note that $h$ (resp.\@ $k$) can possibly be equal to $-\infty$ (resp.\@ $+\infty$)  in which case $\Omega_V$ is a half space.
In particular, $\Omega_V(0,+\infty)$ is the open left half of the complex plane, corresponding to stable matrices for continuous LTI systems.

\item Disks centred on the real line:
the disk centred  at $(-q,0)$ with radius $r>0$, denoted by $\Omega_D(-q,r)$, is defined as
\[
\Omega_D(-q,r)
:= \left\{ z\in \C \ \big| \ |z+q|<r\right\}.
\]
In particular, $\Omega_D(0,1)$ is the unit disk, corresponding to stable matrices for discrete LTI systems.

\end{itemize}

\begin{figure}
\begin{center}
\includegraphics[width=5cm]{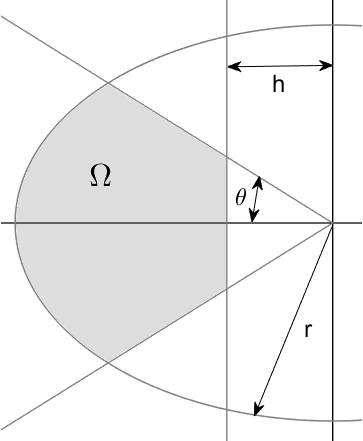}\label{fig1}
\caption{Illustration of $\Omega
 =
\big\{ x+iy \ | \
\sin(\theta) x < \cos(\theta) y < -\sin(\theta) x,\,
x<-h<0,\,
|x+iy|<r
\big\}
=
\Omega_C(0,\theta)
\cap \Omega_V(h,+\infty)
\cap \Omega_D(0,r)$.
}
\end{center}
\end{figure}

The sets $\Omega$ considered in this paper can be either any of $\Omega_C$, $\Omega_V$, $\Omega_D$, or the intersection of such sets; see Figure~\ref{fig1} for an illustration. Note that $\Omega$ is symmetric with respect to the real line. Note also that it is useless to consider more than one set of the type $\Omega_V$ since the intersection of such sets can be simplified to a single set $\Omega_V$.
However, it makes sense to consider the intersection of several sets of the types $\Omega_C$ and $\Omega_D$, that is, the intersection of several conic sectors and several disks.

In this paper, we consider the problem of computing the nearest $\Omega$-stable matrix to a given matrix. More precisely, for a given matrix  $A \in \R^{n,n}$
we are interested in solving the following optimization problem
\begin{equation}\label{eq:prob_def}
\inf_{X\in \mathbb S_{\Omega}^{n,n}} {\|A-X\|}_F^2,
\end{equation}
where ${\|\cdot\|}_F$ denotes the Frobenius norm of a matrix and $\mathbb S_{\Omega}^{n,n}$ is the set of
all $\Omega$-stable matrices of size $n\times n$. This problem is important for example in system identification where
one needs to identify a stable system from observations. For example consider the $\Omega$-stability region as in Figure~\ref{fig1} which is the intersection of the regions $\Omega_C$, $\Omega_V$ and $\Omega_D$.
A solution of the nearest $\Omega$-stable
matrix problem for this region is useful to identify a stable system with a minimum decay rate $h$, a minimum damping ratio
$\tau=\cos\theta$, and a maximum undamped natural frequency $\omega_d=r\sin\theta$. Such a system bounds the
maximum overshoot, the frequency of oscillatory modes, the delay time, the rise time, and the setting time~\cite{ChilG96,Kuo87}.
Such nearness problems for LTI systems have recently been studied; see for example  \cite{OrbNV13,GilS17,MehMS17,GugL17} for the continuous LTI systems, and
\cite{OrbNV13,NesP17,GP2018,GilKS18a} for discrete LTI systems.
The converse of~(\ref{eq:prob_def}) is the distance to $\Omega$-instability, that is, for a given
$\Omega$-stable matrix $A$, find the smallest perturbation $\Delta_A$ with respect to some norm
such that the perturbed matrix $A+\Delta_A$ has at least one eigenvalue outside $\Omega$.
It is the more constrained version of the widely studied distance to instability~\cite{Bye88,HinP86}.

\paragraph{Notation}
Throughout the paper, $X^T$ and
$\|X\|$  stand for the transpose and the spectral norm of a real square matrix $X$, respectively.
We write $X\succ 0$ $(X \preceq 0)$ and $X\succeq 0$ $(X \preceq 0)$ if $X$ is symmetric and positive definite (negative definite)
or positive semidefinite (negative semidefinite), respectively.
By $I_m$ we denote the identity matrix of size $m \times m$.

\paragraph{Outline and contribution}

In order to tackle~\eqref{eq:prob_def}, we extend the approach proposed in~\cite{GilS17} that tackles the nearest stable matrix problem for continuous LTI systems.
In~\cite{GilS17}, stable matrices for continuous LTI systems are parametrized using dissipative Hamiltonian (DH) matrices of the form $A = (J-R)Q$ where $J^T = -J$, $R \succeq 0$ and $Q \succ 0$. In fact, it turns out that a matrix is stable if and only if it is a DH matrix.
In Section~\ref{sec:param}, we show how to impose the eigenvalues of a matrix of the form $A = (J-R)Q$ to lie in the sets $\Omega_C$, $\Omega_V$ and $\Omega_D$
using linear matrix inequalities (LMIs) derived in~\cite{ChilG96}. This allows us to reformulate in Section~\ref{sec:reform} the problem~\eqref{eq:prob_def} into an equivalent optimization problem with a convex feasible set onto which it is easy to project.
We then propose in Section~\ref{sec:algo} a block coordinate descent algorithm to tackle this problem and illustrate the effectiveness of this algorithm on several examples.

\section{DH matrices and $\Omega$-stability} \label{sec:param}

Let us formally define a DH matrix.
\begin{definition}{\rm
A matrix $A \in \R^{n,n}$ is said to be a {\emph DH matrix} if $A=(J-R)Q$ for some $J,R,Q \in \R^{n,n}$
such that $J^T=-J$, $R\succeq 0$ and $Q \succ 0$.}
\end{definition}
It was shown in~\cite{GilS17} that a matrix is stable (that is, all its eigenvalues are in the left half of the complex plane) if and only if it is a DH matrix.
In this section, we obtain a parametrization of the sets $\Omega_C$, $\Omega_V$ and $\Omega_D$ in terms of DH matrices with extra constraint on  $J$, $R$ and $Q$;
see Section~\ref{paramC},~\ref{paramV}
and~\ref{paramD}, respectively.
Note that for $\Omega$, the constraint $R \succeq 0$ can be removed to allow the region to intersect with the right half of the complex plane. This will allow us in particular to model
$\Omega$-stability for discrete time system that correspond to $\Omega_D(0,1)$ or its intersection with $\Omega_C$ and $\Omega_V$.
The following lemma will be frequently used in the following subsections.

\begin{lemma}\label{lem:dhevalues}
Let $A=(J-R)Q$, where $J,R,Q \in \R^{n,n}$ such that $J^T=-J$, $R^T=R$, and $Q^T=Q$ is invertible.
Let $\lambda \in \C$, and $v \in \C^{n}\setminus \{0\}$  be such that
$v^*A=\lambda v^*$. Then
\[
\real{(\lambda)}= -\frac{v^*Rv}{v^*Q^{-1}v} \quad \text{and}\quad \imag{(\lambda)}= -i\frac{v^*Jv}{v^*Q^{-1}v}.
\]
\end{lemma}
\proof
Given $v$ a left eigenvector of $A$ corresponding to eigenvalue $\lambda$, i.e., $v^*A=\lambda v^*$. This implies that
\begin{eqnarray}\label{eq1}
v^*(J-R)Q=\lambda v^* \Longrightarrow v^*(J-R)v=\lambda v^*Q^{-1}v
\end{eqnarray}
and by taking the conjugate of~\eqref{eq1}, we get
\begin{eqnarray}\label{eq2}
 v^*(-J-R)v=\bar \lambda v^*Q^{-1}v.
\end{eqnarray}
From~\eqref{eq1} and~\eqref{eq2} we have
\begin{eqnarray*}\label{eq3}
2v^*Jv &=&(\lambda -\bar \lambda)v^*Q^{-1}v =2 i \lambda_2 v^*Q^{-1}v \Longrightarrow v^*Jv = i \imag{(\lambda)} v^*Q^{-1}v,
\end{eqnarray*}
and
\begin{eqnarray*}\label{eq4}
-2v^*Rv &=&(\lambda +\bar \lambda)v^*Q^{-1}v =2 \lambda_1 v^*Q^{-1}v \Longrightarrow v^*Rv = - \real{(\lambda)} v^*Q^{-1}v.
\end{eqnarray*}
\eproof

\subsection{Parametrization for conic sectors $\Omega_C$} \label{paramC}

Consider the $\Omega_C(a,\theta)$ region with parameters $a\in \R$ and $0 \leq \theta \leq \pi/2$ and
set $\alpha := \sin(\theta)$ and $\beta := \cos(\theta)$.
To parametrize $\Omega_C(a,\theta)$ in terms of DH matrices, let us first prove the following elementary lemma.
\begin{lemma}\label{lem1}
Let $\lambda =\lambda_1+i\lambda_2$, where $\lambda_1$, $\lambda_2 \in \R$. Then  $\mat{cc} \alpha \,(\lambda_1-a) & \beta \,i\lambda_2 \\ -\beta\, i\lambda_2 & \alpha \,(\lambda_1-a) \rix \prec 0$
if and only if $\lambda \in \Omega_C(a,\theta)$.
\end{lemma}
\proof The proof follows using the fact that $\mat{cc} \alpha \,(\lambda_1-a) & \beta \,i\lambda_2 \\ -\beta\, i\lambda_2 & \alpha \,(\lambda_1-a) \rix$
is Hermitian and therefore it is negative definite if and only if both eigenvalues $\mu_1=\alpha (\lambda_1-a) +\beta \lambda_2$
and $\mu_2=\alpha (\lambda_1-a) -\beta \lambda_2$ are negative which is true if and only if $\alpha (\lambda_1-a) < \beta \lambda_2 < -\alpha (\lambda_1-a)$, i.e.,
$\lambda \in \Omega_C(a,\theta)$.
\eproof

\begin{theorem}\label{mainthm}
Let $A\in \R^{n,n}$. Then $A$ is $\Omega_C(a,\theta)$-stable if and only if $A=(J-R)Q$ for some $J,R,Q \in R^{n,n}$ such that
$Q \succ 0$ and
\begin{equation}\label{eq:tem1}
\mat{cc}\alpha  (R+aQ^{-1}) & -\beta J \\
\beta   J & \alpha (R+aQ^{-1})
\rix
\succ 0.
\end{equation}
\end{theorem}
\proof
First suppose that $A=(J-R)Q$ for some $J,R,Q$ satisfying $Q\succ 0$ and~\eqref{eq:tem1}.
Let $\lambda =\lambda_1 +i \lambda_2$ be an eigenvalue of $A$ and let $v \in \C^n \setminus \{0\}$ be a left eigenvector of $A$ corresponding to eigenvalue $\lambda$.
Since $\mat{cc}\alpha\, (R+aQ^{-1}) & -\beta\,J \\\beta\, J & \alpha\, (R+aQ^{-1})\rix \succ 0$ and $v \neq 0$, we have that
\begin{eqnarray}\label{eq51}
-2\mat{cc} v^* & 0 \\ 0 & v^*\rix
\mat{cc}\alpha\, (R+aQ^{-1}) & -\beta\,J \\\beta\, J & \alpha\, (R+aQ^{-1})\rix
\mat{cc} v & 0 \\ 0 & v\rix &\prec& 0 \nonumber\\
\Longrightarrow\quad  2\mat{cc} -\alpha v^*(R+aQ^{-1})v & \beta v^*Jv \\-\beta v^*Jv & -\alpha v^*(R+aQ^{-1})v \rix &\prec& 0 \nonumber\\
\Longrightarrow\quad  2\mat{cc} -\alpha v^*Rv  & \beta v^*Jv \\-\beta v^*Jv & -\alpha v^*Rv \rix -
\alpha \mat{cc}av^*Q^{-1}v & 0 \\ 0 & av^*Q^{-1}v \rix&\prec& 0.
\end{eqnarray}
Thus by using Lemma~\ref{lem:dhevalues} in~\eqref{eq51}, we obtain
\begin{eqnarray}
v^*Q^{-1}v\mat{cc} \alpha \,(\lambda_1-a) & \beta\,i\lambda_2 \\-\beta\,i \lambda_2 & \alpha \,(\lambda_1-a) \rix \prec 0.
\end{eqnarray}
This implies that $\mat{cc} \alpha \,(\lambda_1-a) & \beta\,i\lambda_2 \\-\beta\,i \lambda_2 & \alpha \,(\lambda_1-a) \rix \prec 0$ since $Q$ is positive definite.
Thus Lemma~\ref{lem1} implies that $A$ is $\Omega_C(a,\theta)$-stable.

For the `only if' part, since $A$ is $\Omega_C(a,\theta)$-stable, by~\cite[Theorem~2.2]{ChilG96}, there exists $X \succ 0$ such that
\begin{equation}\label{eqt}
\mat{cc} \alpha (AX+XA^T -aX) & \beta (AX-XA^T) \\
\beta (XA^T-AX) & \alpha (AX+XA^T-aX)
\rix \prec 0.
\end{equation}
Let
\[
R=-\frac{AX+(AX)^T}{2},\quad J=\frac{AX-(AX)^T}{2},\quad \text{and}\quad Q=X^{-1}.
\]
Then $(J-R)Q=A$ and it follows from~\eqref{eqt} that
\[
\mat{cc}\alpha\, (R+aQ^{-1}) & -\beta\,J \\\beta\, J & \alpha\, (R+aQ^{-1})\rix =- 1/2 \mat{cc} \alpha (AX+XA^T-aX) & \beta (AX-XA^T) \\
\beta (XA^T-AX) & \alpha (AX+XA^T-aX) \rix \succ 0.
\]
\eproof

As a consequence of~\eqref{eq:tem1} in Theorem~\ref{mainthm}, the matrix $J$ is skew-symmetric. However, the matrix $R$ may not be positive definite (when $a >0$)
and therefore the $\Omega_C(a,\theta)$-stable matrix $A$ need not be a DH matrix. But
when $a \leq 0$, then~\eqref{eq:tem1} implies that $R+aQ^{-1} \succ 0$, or equivalently, $R \succ -aQ^{-1}$ since $a \leq 0$ and $Q\succ 0$. As a result
$R$ is positive semidefinite. Therefore in this case $A$ is $\Omega_C(a,\theta)$-stable if and only if $A$ is a DH matrix satisfying~\eqref{eq:tem1}.

\subsection{Parametrization for vertical strips $\Omega_V$} \label{paramV}

We can characterize $\Omega_V$-stability as follows.

\begin{theorem}\label{thm:vert}
Let $A \in \R^{n,n}$ and $ h < k$. Then $A$ is $\Omega_V(h,k)$-stable if and only if $A=(J-R)Q$
for some $J,R,Q \in \R^{n,n}$ such that $J^T=-J$, $R^T=R$, $Q\succ 0$, and
\begin{equation} \label{condVert}
k Q^{-1} \succ R \succ h Q^{-1}.
\end{equation}
\end{theorem}
\proof First suppose that $A=(J-R)Q$, where $J^T=-J$, $R^T=R$, $Q\succ 0$ such that  $k Q^{-1} \succ R\succ h Q^{-1}$.
Let $\lambda$ be an eigenvalue of $A$ and $x \in \C^{n}\setminus \{0\}$ such that
$Ax=\lambda x$ or $(J-R)Qx=\lambda x$. Since $Q$ is invertible, this implies that
\begin{equation}\label{eq:1}
x^*Q(J-R)Qx=\lambda x^*Q x \quad \Longrightarrow \quad \real{(\lambda)}=-\frac{x^*QRQx}{x^*Qx}.
\end{equation}
Since $x^*Qx > 0$ as $Q\succ 0$ and $R$ satisfies $k Q^{-1} \succ R\succ h Q^{-1}$, we have
$k x^*Qx > x^*QRQx > h x^*Qx$. This implies that
\begin{equation}\label{eq:2}
k  > \frac{x^*QRQx}{x^*Qx} > h.
\end{equation}
From~\eqref{eq:1} and~\eqref{eq:2}, we have that $-k<\real{(\lambda)} <-h$.

Conversely, let $A$ be $\Omega_V(h,k)$-stable. Then from~\cite{ChilG96}, there exists a single $X \succ 0$
such that
\begin{equation}\label{eq:3}
AX+XA^T+2hX \prec 0\quad  \text{and}\quad AX+XA^T+2kX \succ 0.
\end{equation}
Define
\[
Q=X^{-1},\quad R=-\frac{AX+XA^T}{2},\quad J=\frac{AX-XA^T}{2}.
\]
Then clearly $A=(J-R)Q$. Also in view of~\eqref{eq:3} $R$ satisfies $k Q^{-1} \succ R\succ h Q^{-1}$.
\eproof

It is easy to see that in Theorem~\ref{thm:vert} when $h \geq 0$ $A$ is $\Omega_V$-stable
if and only if $A$ is a DH matrix since $R \succ h Q^{-1}$.

\subsection{Parametrization for disks $\Omega_D$}  \label{paramD}

We first recall a result from~\cite{ChilG96} and state it as a lemma in the following.
\begin{lemma}{\rm \cite{ChilG96}}\label{lem:omega_stab}
Consider the $\Omega_D(-q,r)$ region where $q \in \R$ and $r > 0$, and let $\lambda \in \C$.
Then  $\lambda \in \Omega_D(-q,r)$ if and only if
$\mat{cc} -r & q+\lambda \\ q+\overline{\lambda} & -r \rix \prec 0$.
\end{lemma}
We can characterize $\Omega_D$-stability as follows.
\begin{theorem}\label{thm:omega2}
Let $A \in \R^{n,n}$, $q \in \R$ and $r > 0$.
Then $A$ is $\Omega_D(-q,r)$-stable if and only if  $A=(J-R)Q$ for some $J,R,Q \in \R^{n,n}$ such that
$J^T = -J$, $R$ is symmetric, $Q \succ 0$ and
\begin{equation}\label{eq:t1}
\mat{cc} rQ^{-1}& -qQ^{-1} \\ -qQ^{-1} & r Q^{-1} \rix \succ \mat{cc} 0& J-R \\ (J-R)^T &0
\rix.
\end{equation}
\end{theorem}
\proof 
First suppose that $A=(J-R)Q$ with $J^T=-J$, $R^T=R$ and $Q\succ 0$ such that~\eqref{eq:t1} holds.
Let $\lambda=\lambda_1+i\lambda_2$, where $\lambda_1,\lambda_2 \in \R$ be an eigenvalue of $A$ and $v \in \C^{n}\setminus \{0\}$ be a corresponding left eigenvector.
Since~\eqref{eq:t1} holds, we have that
\begin{eqnarray*}
\mat{cc} v^* & 0 \\ 0 & v^*\rix\mat{cc} rQ^{-1}& -qQ^{-1} \\ -qQ^{-1} & r Q^{-1} \rix
\mat{cc} v & 0 \\ 0 & v\rix \succ
\mat{cc} v^* & 0 \\ 0 & v^*\rix\mat{cc} 0& J-R \\ (J-R)^T &0
\rix \mat{cc} v & 0 \\ 0 & v\rix.
\end{eqnarray*}
This implies that
\begin{eqnarray*}
\mat{cc} r & -q \\ -q & r  \rix v^*Q^{-1}v \succ \mat{cc} 0& v^*(J-R)v \\ v^*(J-R)^Tv &0
\rix.
\end{eqnarray*}
Since $v^*Q^{-1}v > 0$ as $Q \succ 0$, we obtain
\begin{eqnarray*}
 \mat{cc} r & -q \\ -q & r  \rix  \succ \mat{cc} 0& \frac{v^*(J-R)v}{v^*Q^{-1}v} \\ \frac{v^*(J-R)^Tv}{v^*Q^{-1}v} &0
\rix .
\end{eqnarray*}
Thus in view of Lemma~\ref{lem:dhevalues}, we have
\begin{eqnarray*}
 \mat{cc} r & -q \\ -q & r  \rix \succ \mat{cc} 0& i\lambda_2+\lambda_1 \\ -i\lambda_2+\lambda_1 &0
\rix \quad \Longrightarrow \quad  \mat{cc} r & -q-\lambda \\ -q-\overline{\lambda} & r  \rix \succ 0.
\end{eqnarray*}
This implies by using Lemma~\ref{lem:omega_stab} that $\lambda \in \Omega_D(-q,r)$ and therefore $A$ is
$\Omega_D(-q,r)$-stable.

Conversely, suppose $A$ is $\Omega_D(-q,r)$-stable then by \cite[Theorem~2.2]{ChilG96} there exists $X \succ 0$ satisfying
\begin{equation}\label{eqtr1}
\mat{cc} -rX & qX+AX \\ qX +XA^T & -rX
\rix \prec 0.
\end{equation}
Define
\[
Q=X^{-1},\quad R=-\frac{AX+XA^T}{2},\quad J=\frac{AX-XA^T}{2}.
\]
Then $A=(J-R)Q$ with $J^T = -J$, $R$ symmetric and $Q \succ 0$.
Moreover, by~\eqref{eqtr1}, we have
\begin{eqnarray*}
0 \succ \mat{cc} -rX & qX+AX \\ qX +XA^T & -rX \rix &=& \mat{cc} -rQ^{-1} & qQ^{-1}+AQ^{-1} \\ qQ^{-1} +Q^{-1}A^T & -rQ^{-1}
\rix \\
&=&\mat{cc} -r Q^{-1}  & q Q^{-1}  \\ q Q^{-1}  & -r Q^{-1}  \rix  +
\mat{cc}0 & J-R \\ (J-R)^T & 0\rix.
\end{eqnarray*}
This implies that
\begin{eqnarray*}
\mat{cc} r Q^{-1}  & -q Q^{-1}  \\ -q Q^{-1} & r Q^{-1}  \rix \succ
\mat{cc}0 & J-R \\ (J-R)^T & 0\rix.
\end{eqnarray*}
This completes the proof.
\eproof

We note that in the above theorem, the matrix $R$ need not be positive semidefinite and thus a $\Omega_D$-stable matrix need not be a DH matrix. However, if the disc $\Omega_D$ completely lies in the left half of the complex plane,
then $A$ is a DH matrix. More precisely we have the following result.
\begin{theorem}
Let $A \in \R^{n,n}$ and $q \geq r > 0$.
Then $A$
is $\Omega_D(-q,r)$-stable if and only if $A$ is a DH matrix such that~\eqref{eq:t1} is satisfied.
\end{theorem}
\proof The proof is immediate from Theorem~\ref{thm:omega2} if we show that $A$ is $\Omega_D(-q,r)$-stable implies $R=-\frac{AX+XA^T}{2}$ is positive semidefinite. Let $x\in \C^n\setminus\{0\}$ so that $2x^*Rx=-x^*(AX)x-x^*(AX)^Tx$.
Using~\eqref{eqtr1}, we also have
\begin{eqnarray*}
&\mat{c} x \\ x\rix^* \mat{cc} rX & -qX-AX \\ -qX -XA^T & rX \rix \mat{c}x \\x \rix >0.
\end{eqnarray*}
This implies that
$-x^*(AX)x-x^*(AX)^Tx ~>~ 2(q-r)x^*Xx $ and hence $ x^*Rx > 2(q-r)x^*Xx \geq 0,$
since $X \succ 0$ and $q \geq r$ which completes the proof.
\eproof

\section{Reformulation of the nearest $\Omega$-stable matrix problem} \label{sec:reform}

In this section, we reformulate the problem~\eqref{eq:prob_def} of finding the nearest $\Omega$-stable matrix to a given unstable matrix $A$ into an equivalent optimization problem with a convex feasible set.
In the reformulations of $\Omega_V$ and $\Omega_D$, we used linear matrix inequalities (LMIs) involving the inverse of $Q^{-1}$, which is not convex; see~\eqref{eq:tem1},~\eqref{condVert} and~\eqref{eq:t1}.
Let us introduce the auxiliary variable $P = Q^{-1} \succ 0$.
In view of Theorems~\ref{mainthm},~\ref{thm:vert} and~\ref{thm:omega2}, this allows us to parametrize the sets
$\mathbb S_{\Omega_C(a,\theta)}^{n,n}$,
$\mathbb S_{\Omega_V(h,k)}^{n,n}$, and
$\mathbb S_{\Omega_D(-q,r)}^{n,n}$ as follows:
\begin{equation} \label{eq:omegaC}
\mathbb S_{\Omega_C(a,\theta)}^{n,n} = \left\{(J-R)P^{-1} \ \big| \ J,R,P \in \R^{n,n}, P \succ 0,~\mat{cc}
\sin(\theta) \, (R+aP) & -\cos(\theta) \,J \\
\cos(\theta) \, J & \sin(\theta) \, (R+aP)
\rix \succ 0
\right\},
\end{equation}
\begin{equation} \label{eq:omegaV}
\mathbb S_{\Omega_V(h,k)}^{n,n} = \left\{(J-R)P^{-1} \ \big| ~ J,\,R,\,P \in \R^{n,n},\, P \succ 0,\, k P \succ R\succ h P
\right\},
\end{equation}
and
\begin{equation} \label{eq:omegaD}
\mathbb S_{\Omega_D(-q,r)}^{n,n} = \left\{(J-R)P^{-1} \ \big| \ J,R,P \in \R^{n,n}, J^T = -J,
\mat{cc} rP & -qP-(J-R) \\ -qP-(J-R)^T & r P \rix \succ 0
\right\},
\end{equation}
where $0 \leq \theta < \pi/2$, $ h < k$ and $r > 0$.
Note that these sets are non-convex and open.  For example, the set of stable matrices (or, equivalently $\Omega_C(0,\theta)$ matrices with $\theta=\pi/2$) is non-convex~\cite{OrbNV13}.
It can be checked that for any $\theta \in (0,\pi/2)$, the set $\Omega_C(0,\theta)$ is also non-convex. Indeed, consider $A=\mat{cc} \alpha \,(-\beta) & \beta \,i (\alpha) \\ 0 & \alpha \,(-\beta) \rix $
and $B=\mat{cc} \alpha \,(-\beta) & 0\\ -\beta \,i (\alpha) & \alpha \,(-\beta) \rix $ with $A,\,B \in \mathbb S_{\Omega_\theta}^{n,n}$. For $\gamma=1/2$, the matrix $\gamma A+(1-\gamma)B \notin \mathbb S_{\Omega_\theta}^{n,n}$ as it has an eigenvalue at zero.
Note also that the above sets are not closed (because of the constraints of positive definiteness). From an optimization point of view, it does not make much sense to optimize on such sets since the optimal solution(s) may not be attained. Therefore, we will consider the closure of these sets: this amounts to replace all constraints involving a positive definite constraint with a positive semidefinite  constraint, that is,
replace $\succ 0$ with $\succeq 0$,
in the definition of the sets~\eqref{eq:omegaC}, \eqref{eq:omegaV} and~\eqref{eq:omegaD}.
We will denote the corresponding sets as
 $\bar{\mathbb S}_{\Omega_C(a,\theta)}^{n,n}$,
$\bar{\mathbb S}_{\Omega_V(h,k)}^{n,n}$, and
$\bar{\mathbb S}_{\Omega_D(-q,r)}^{n,n}$, respectively. Note that by considering the closure of these sets, as done in~\cite{GilS17}, we do no change the value of the infimum of~\eqref{eq:prob_def}.

Finally, given $0 < \theta \leq \pi/2$, $ h < k$ and several disks of parameters $(q_i,r_i)$ $1 \leq i \leq k$,  we tackle~\eqref{eq:prob_def} by solving
\begin{equation} \label{finalform}
\inf_{J,R,P} {\|A - (J-R)P^{-1}\|}_F^2
\quad
\text{ such that }
\quad
 (J-R)P^{-1}
\in  \bar{\mathbb S}_\Omega^{n \times n},
\end{equation}
where
\[
\bar{\mathbb S}_\Omega^{n \times n}
\; = \;
\bar{\mathbb S}_{\Omega_C(a,\theta)}^{n,n}
\cap
\bar{\mathbb S}_{\Omega_V(h,k)}^{n,n}
\cap_{i=1}^k
\bar{\mathbb S}_{\Omega_D(q_i,r_i)}^{n,n}.
 \]
 The feasible set of the above optimization problem only involves convex LMI constraints. Of course, the objective function is non-convex and the problem remains difficult,
 but it is easier to handle non-convex objective rather than a non-convex feasible set.

\section{Algorithm and numerical experiments} \label{sec:algo}

In this section, we describe a block coordinate descent method to tackle~\eqref{finalform} and then apply it on some selected examples.
To solve the convex optimization subproblems involving LMIs, we use the interior point method SDPT3 (version 4.0)~\cite{toh1999sdpt3, tutuncu2003solving} with CVX as a modeling system~\cite{cvx, gb08}.
Our code is available from \url{https://sites.google.com/site/nicolasgillis/} and the numerical examples presented below can be directly run from this online code.
All tests are preformed using Matlab
R2015a on a laptop Intel CORE i7-7500U CPU @2.7GHz 24Go RAM.

\subsection{Block coordinate descent method for~\eqref{finalform}}

In this section, we propose to use a block coordinate descent (BCD) method to tackle~\eqref{finalform}.
The BCD method will alternatively optimize the block of variable $(J,R)$ and $P$:
\begin{itemize}

\item For $P$ fixed, the subproblem in $(J,R)$ is convex since ${\|A-(J-R)P^{-1}\|}_F^2$
is a quadratic function of $(J,R)$, and the feasible set is convex in variables $(J,R,P)$. More precisely, this is a linear least squares problem under LMI constraints.

\item For $(J,R)$ fixed, the problem in variable $P$ is non-convex, because of the objective function ${\|A-(J-R)P^{-1}\|}_F^2$. In a standard BCD method, one would optimize solely on variable $P$ for $(J,R)$ fixed. However, we have observed that performing a gradient step on all the variables performs better. Denoting $f(J,R,P) = {\|A - (J-R)P^{-1}\|}_F^2$ and
$E = (J-R)P^{-1}-A$,
the gradient of $f$ is given by
\begin{align*}
\nabla_J f(J,R,P) & = - \nabla_R f(J,R,P) = 2 E P^{-T},  \\
\nabla_P f(J,R,P) & = - 2 P^{-T} (J-R)^T E P^{-T}.
\end{align*}
For the computation of the gradient with respect to $P$,
we refer to~\cite[Appendix]{GilMS17} for a similar derivation.
For the step length used in the gradient method,
we use a backtracking line search, that is, we reduce the step length as long as the objective function increases, while allowing it to increase at the next step.

\end{itemize}

Note that we have also implemented a projected fast gradient method, as done in~\cite{GilS17}, but it does not perform as well as BCD. The reason is that, in~\cite{GilS17}, the feasible set is simply $R \succeq 0$ and $Q=P^{-1}\succeq 0$ which can be projected onto very efficiently with an eigenvalue decomposition.
For general $\Omega$ stability, projecting $(J,R,P)$ onto the feasible set roughly has the same computational cost as optimizing exactly over variables $(J,R)$ for $P$ fixed. 
For all the numerical experiments presented below, we will use 100 iterations of the above scheme.

\subsection{Initialization} \label{init}

As expected, our algorithm will be sensitive to the initial choice of the matrices $(J,R,P)$ since we are tackling the difficult non-convex optimization problem~\eqref{finalform}.
In this paper, we propose two initializations.

The first initialization, which we will refer to as identity initialization chooses $P = I_n$ and then, for $P$ fixed, chooses $J$ and $R$ that minimize~\eqref{finalform}. We will denote this initialization $(J_i,R_i,P_i)$ and the solution obtained by BCD with this initialization as $(J_i^*,R_i^*,P_i^*)$.

The second initialization, which we will refer to as the LMI initialization, will use a solution to a relaxation of the LMIs~\eqref{eqt}, \eqref{eq:3} and \eqref{eqtr1}.
In fact, since the input matrix $A$ is in general not $\Omega$-stable, these LMIs do not admit a feasible solution $X \succ 0$. Hence we will replace all inequalities $\succ 0$ with $\succ -\delta I$ and minimize $\delta$ over variables $X \succ I$ and $\delta$. Denoting $(X^*,\delta^*)$ an optimal solution of these relaxed LMIs, we take
\[
P=X^*,\quad R=-\frac{AX^*+X^*A^T}{2},\quad J=\frac{AX^*-X^*A^T}{2},
\]
as specified in
Theorems~\ref{mainthm},~\ref{thm:vert} and~\ref{thm:omega2}.
Note that the optimal solutions $(X^*,\delta^*)$ will be such that $\delta^* = 0$ if and only if $A$ is $\Omega$-stable in which case the problem would be solved.
This is an advantage of this initialization: it will detect whether the input matrix is $\Omega$-stable.
We will denote this initialization $(J_x,R_x,P_x)$ and the solution obtained by BCD with this initialization as $(J_x^*,R_x^*,P_x^*)$.

\subsection{Synthetic data sets}

\begin{figure}
\begin{center}
\begin{tabular}{cc}
\includegraphics[height=0.42\textwidth]{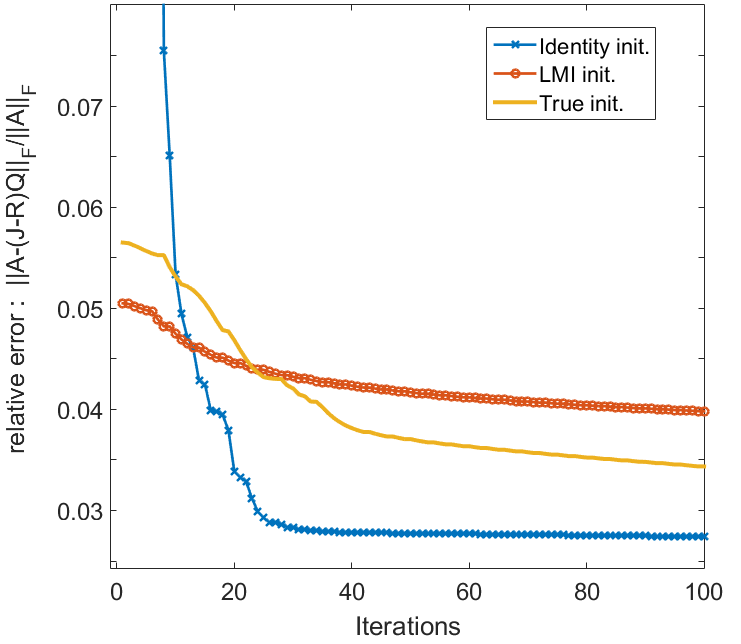}
&
\includegraphics[height=0.42\textwidth]{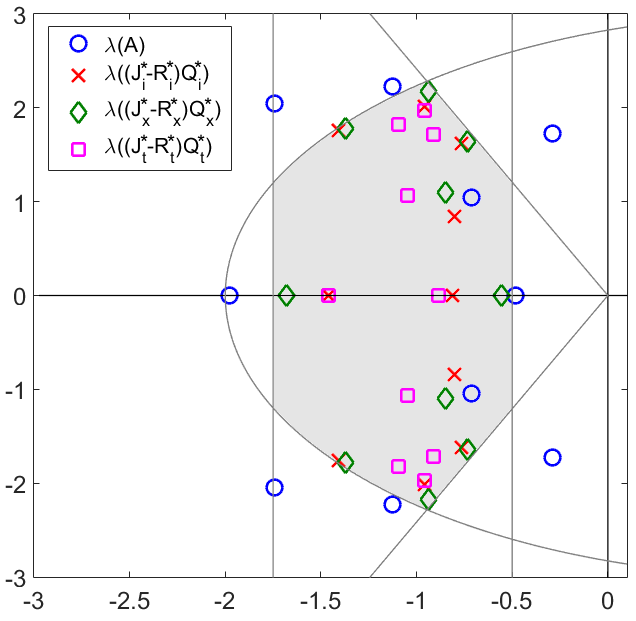}
\end{tabular}
\caption{ On the left: evolution of the relative error for the different initializations.
On the right: eigenvalues of $A$, $(J_i^*-R_i^*)Q_i^*$, $(J_x^*-R_x^*)Q_x^*$ and $(J_t^*-R_t^*)Q_t^*$.
The shaded area is the set
$\Omega =
\Omega_C(3\pi/8)
\cap \Omega_V(0.5,1.75)
\cap \Omega_D(1,3)$. \label{fig:synth} }
\end{center}
\end{figure}

Let us consider
\[
\Omega =
\Omega_C(0,3\pi/8)
\cap \Omega_V(0.5,1.75)
\cap \Omega_D(-1,3),
\]
so that
\[
\theta = \frac{3\pi}{8},
h = 0.5,
k = 1.75,
q = -1 \text{ and } r = 3.
\]
The set $\Omega$ is illustrated on Figure~\ref{fig:synth}.
To generate a matrix that is not too far from being $\Omega$-stable, we proceed as follows.
First, we generate each entry of $J_0$, $R_0$ and $Q_0$ at random using the Gaussian distribution of mean 0 and standard deviation 1 (\texttt{randn(n)} in Matlab). Then we perform the following projection
\[
(J_t,R_t,P_t)
\quad
=
\quad
\argmin_{J,R,P} {\|(J,R,P)-(J_0,R_0,P_0)\|}_F \text{ such that } (J-R)P^{-1} \in \bar{\mathbb S}_\Omega^{n \times n},
\]
to obtain an $\Omega$-stable matrix $A_t = (J_t-R_t)P_t^{-1}$.
Finally, we generate a matrix $N$ where each entry is generated at random using the Gaussian distribution of mean 0 and standard deviation 1, and set
\[
{A} = A_t + \epsilon {\|A_t\|}_F \frac{N}{{\|N\|}_F},
\]
so that ${\|A_t-{A}\|}_F/{\|A_t\|}_F = \epsilon$. The matrix $A$ will in general not be $\Omega$-stable for $\epsilon$ sufficiently large.

In the following, we compare the solutions obtained by our BCD algorithm with the two initializations described in the previous section. We also consider the initialization made of the true $(J_t,R_t,Q_t)$ used to generate $A$, which we will refer to as the true initialization. This will be useful to validate the other two initializations and see whether they are able to obtain good nearby $\Omega$-stable matrices.

Figure~\ref{fig:synth} displays the result of a particular example with $\epsilon = 0.1$. We observe that the three initializations lead to different solutions. The best solution (with relative error $\frac{{\|A-(J-R)P^{-1}\|}_F}{{\|A\|}_F} = 2.74\%$) is found with the identity initialization although it has, as expected, the highest initial relative error (41.33\%). The LMI initialization has the lowest initial error (5.04\%), lower than the true initialization\footnote{The true initialization
$(J_t-R_t)P_t^{-1} = A$ has relative error $\epsilon=10\%$ by construction. However, in BCD, we only use $P$ for the initialization and set $(J,R)$ as the optimal solution for $P$ fixed.} (5.64\%) but BCD is not able to improve it much (up to 3.9\%).


Let us now perform more extensive numerical experiments with these randomly generated matrices.
We consider $\epsilon = 0.01, 0.05, 0.1, 0.2$, and for each value of $\epsilon$, we generate 10 matrices $A$.
Table~\ref{tab1} reports the average error and standard deviation for each noise level for these 10 randomly generated matrices. It also reports in bracket the number of times each initialization obtained the best solution (up to 0.01\%).
\begin{center}
 \begin{table}[h!]
 \begin{center}
\caption{Comparison of the algorithms for randomly generated matrices $A$.
The table displays  the average relative error in percent with the standard deviation obtained by each initialization and, in brackets,
the number of times the algorithm found the best solution out of the 20 runs (up to 0.01\%). The best results are highlighted in bold. \label{tab1}}
 \begin{tabular}{|c|c|c|c|}
 \hline   &  Identity & LMI & True    \\
 \hline \hline
 $\epsilon = 0.01$
 &  0.53 $\pm$  0.70 (18) &  \textbf{0.00} $\pm$ 0.01  (\textbf{20}) &  0.21 $\pm$  0.12 (0)\\ \hline
 $\epsilon = 0.05$ &  0.99 $\pm$  0.53 (14) &  \textbf{0.25} $\pm$  0.40 (\textbf{19}) &  1.37 $\pm$  0.74 (1)\\ \hline
 $\epsilon = 0.10$ &  \textbf{2.06} $\pm$  1.64 (\textbf{18}) &  2.09 $\pm$  2.24 (12) &  2.76 $\pm$  1.63 (0)\\ \hline
 $\epsilon = 0.20$ &  4.93 $\pm$  2.01 (\textbf{16}) &  \textbf{4.79} $\pm$  2.35 (11) &  6.29 $\pm$  2.63 (0)\\ \hline
\end{tabular}
 \end{center}
 \end{table}
 \end{center}

We observe the following
\begin{itemize}

\item For low noise levels ($\epsilon \leq 0.05$), the LMI initialization performs in average the best, although the identity initialization finds in most of the cases the best solution. This means that, in several cases, the identity initialization is not able to find a good solution.
Hence the LMI initialization is more reliable (the standard deviation is smaller for $\epsilon \leq 0.05$).
This was expected as the LMI initialization provides an exact solution for $\Omega$-stable matrices and, for small noise levels, it is expected that the relaxed LMIs are close to the original LMIs.
Also, recall that the LMI initialization obtains a low error prior to BCD performing any iteration while the identity initialization requires several iterations before achieving a low error;
see Figure~\ref{fig:synth} for an illustration.

\item For larger noise levels ($\epsilon \geq 0.1$), the identity and LMI initializations perform similarly, although the identity initialization finds in more cases the best solution.

\item In all cases, rather surprisingly, the true initialization does not perform well.
For $\epsilon = 0.01$, it performs better than the identity initialization while, in all other cases, it is worse than the other two initializations.
This is rather surprising, and shows that our two proposed initializations are performing well.


\end{itemize}

\subsection{Discrete-time stability: $\Omega = \Omega_D(0,1)$}

\begin{figure}
\begin{center}
\includegraphics[width=0.7\textwidth]{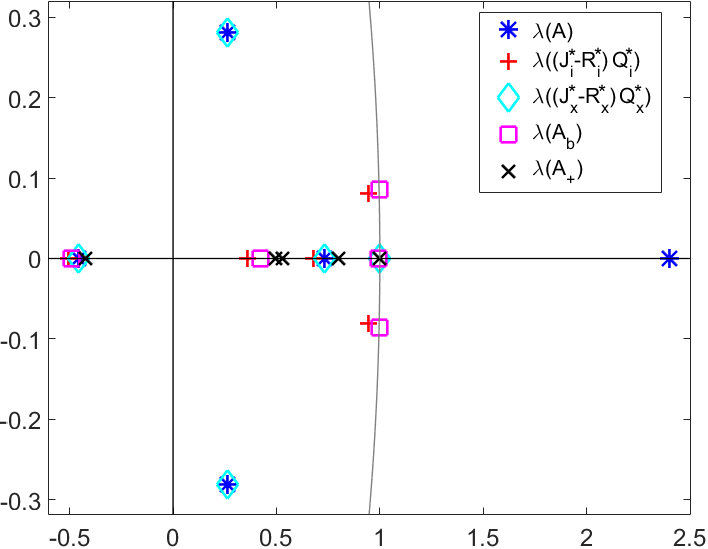}
\caption{ Eigenvalues of $A$, $A_+$, $A_b$, $(J_i^*-R_i^*)Q_i^*$ and $(J_x^*-R_x^*)Q_x^*$. \label{fig:gugl}}
\end{center}
\end{figure}

Stability of discrete-time LTI systems requires the eigenvalues of the matrix $A$ defining the system to belong to the unit disk, that is, to $\Omega_D(0,1)$; see, e.g.,~\cite{OrbNV13, GP2018, GilKS18a} and the references therein.
Therefore, our algorithm can be used to find a nearby stable system for discrete LTI systems. We illustrate this on the example from~\cite[Section 4.4]{GP2018}:
\[
A =
\left( \begin{array}{ccccc}
 0.7 &  0.2 &  0.1 &  0.5 &  1 \\
 0.3 &  0.6 &  0.2 &  0.8 &  0.3 \\
 0.5 &  0.7 &  0.9 &  1 &  0.5 \\
 0.1 &  0.1 &  0.3 &  0.8 &  0.3\\
 0.8 &  0.2 &  0.9 &  0.3 &  0.2 \\
\end{array} \right)
\]
\normalsize
with $\rho(A) = \max_i |\lambda_i(A)| = 2.4$. The nonnegative solution provided by the authors with their algorithm is
\[
A_+ =
\left( \begin{array}{ccccc}
 0.3796 &  0.1797 &  0 &  0.5 &  0.7343 \\
 0 &  0.5791 &  0.0069 &  0.8 &  0.0274 \\
 0.0580 &  0.6719 &  0.6403 &  1 &  0.1334 \\
 0 &  0 &  0 &  0.8 &  0 \\
 0.4204 &  0.1759 &  0.6770 &  0.3 &  0 \\
\end{array} \right)
\]
\normalsize with error ${\|A-A_+\|}_F = 1.10$.
In~\cite{GilKS18a}, the best reported solution is
\[
A_b =
\left( \begin{array}{ccccc}
0.5999 & 0.1317 &-0.0882& 0.5337& 0.8834\\
    0.2582& 0.5864 &0.0967 &0.8512 &0.2089\\
    0.4469&  0.6904 &0.8242 &1.0419& 0.4257\\
    -0.0828  &-0.1243 &-0.2132 &0.8209& 0.0595\\
    0.7076 &0.1273 & 0.7126 &0.3255 &0.0923
\end{array} \right)
\]
with error ${\|A-A_b\|}_F =0.76$. The solution provided by BCD with the identity initialization is
\[
(J_i^*-R_i^*) Q_i^* =
\left( \begin{array}{ccccc}
    0.5736 &   0.1019 &  -0.1508 &   0.4787 &   0.8503\\
    0.2450 &   0.5561 &   0.0914 &   0.7978 &   0.2382\\
    0.3846 &   0.5947 &   0.6766 &   1.0656 &   0.4030\\
   -0.1054 &  -0.0284  & -0.1183 &   0.5988 &  -0.0210\\
    0.6210 &   0.0515 &   0.5482 &   0.3217 &   0.0121
   \end{array} \right)
\]
    with error ${\|A-(J-R)Q\|}_F =0.90$,
while with the LMI initialization, BCD gives
\[
(J_x^*-R_x^*) Q_x^* =
\left( \begin{array}{ccccc}
    0.4643  &  0.0320 &  -0.1498 &   0.1425  &  0.7628\\
    0.1310  &  0.4773 &   0.0193 &   0.5412  &  0.1311\\
    0.2496  &  0.5191 &   0.6329 &   0.6176  &  0.2492\\
   -0.2572   &-0.1579 &  -0.0810 &   0.2548  & -0.0578\\
    0.5650  &  0.0300 &   0.6486 &  -0.0573  & -0.0350
   \end{array} \right)
\]
with error ${\|A-(J_x^*-R_x^*)Q_x^*\|}_F = 1.40$. Figure~\ref{fig:gugl} displays the eigenvalues of the different solutions: again we see that depending on the initialization and the algorithm used, we obtained rather different solutions.
Surprisingly, although
the LMI initialization achieves the highest approximation error of $A$, it provides an optimal approximation of its eigenvalues: four are perfectly recovered while the last one is approximated by its projection on the unit disk. Hence, although this locally optimal solution does not have a low Frobenius norm error, it has an interesting structure.

\section{Conclusion and further work}

In this paper, we have proposed a new parametrization of $\Omega$-stable matrices, where $\Omega$ is the intersection of several regions of the complex plane, namely conic sectors, vertical strips and disks centred on the real line. This allowed us to propose an algorithm to tackle the nearest $\Omega$-stable matrix problem where one is given a matrix $A$ and is looking for the nearest $\Omega$-stable matrix. We illustrated the effectiveness of this approach on several examples.

Further work include the design of faster algorithms to tackle~\eqref{finalform}. In fact, our BCD algorithm  currently  relies on an interior point method to solve the LMIs which does not scale well. For example, on a standard laptop, $n$ can be up to about 50 for which one iterations takes about 1 minute.
Another direction of research is to extended our approach to other regions of the complex plane. For example, LMI regions, that is, subsets of the complex plane that is representable by LMIs~\cite{ChilG96}, would be of particular interest.

\small

\bibliographystyle{siam}
\bibliography{ChouGS}

\end{document}